\documentclass{amsart}

\usepackage{amsfonts,amssymb,verbatim,amsmath,amsthm,latexsym,textcomp,amscd}
\usepackage{latexsym,amsfonts,amssymb,epsfig,verbatim}
\usepackage{amsmath,amsthm,amssymb,latexsym,graphics,textcomp}
\usepackage{graphicx}
\usepackage{color}
\usepackage{url}
\input{psfig.sty}

\begin{document}

\newtheorem{theorem}{Theorem}[section]
\newtheorem{prop}[theorem]{Proposition}
\newtheorem{lemma}[theorem]{Lemma}
\newtheorem{cor}[theorem]{Corollary}
\newtheorem{defn}[theorem]{Definition}
\newtheorem{conj}[theorem]{Conjecture}
\newtheorem{claim}[theorem]{Claim}

\newcommand{\boundary}{\partial}
\newcommand{\bbC}{{\mathbb C}}
\newcommand{\bbD}{{\mathbb D}}
\newcommand{\bbH}{{\mathbb H}}
\newcommand{\bbZ}{{\mathbb Z}}
\newcommand{\bbN}{{\mathbb N}}
\newcommand{\bbQ}{{\mathbb Q}}
\newcommand{\bbR}{{\mathbb R}}
\newcommand{\proj}{{\mathbb P}}
\newcommand{\lhp}{{\mathbb L}}
\newcommand{\tube}{{\mathbb T}}
\newcommand{\cusp}{{\mathbb P}}
\newcommand\AAA{{\mathcal A}}
\newcommand\BB{{\mathcal B}}
\newcommand\CC{{\mathcal C}}
\newcommand\DD{{\mathcal D}}
\newcommand\EE{{\mathcal E}}
\newcommand\FF{{\mathcal F}}
\newcommand\GG{{\mathcal G}}
\newcommand\HH{{\mathcal H}}
\newcommand\II{{\mathcal I}}
\newcommand\JJ{{\mathcal J}}
\newcommand\KK{{\mathcal K}}
\newcommand\LL{{\mathcal L}}
\newcommand\MM{{\mathcal M}}
\newcommand\NN{{\mathcal N}}
\newcommand\OO{{\mathcal O}}
\newcommand\PP{{\mathcal P}}
\newcommand\QQ{{\mathcal Q}}
\newcommand\RR{{\mathcal R}}
\newcommand\SSS{{\mathcal S}}
\newcommand\TT{{\mathcal T}}
\newcommand\UU{{\mathcal U}}
\newcommand\VV{{\mathcal V}}
\newcommand\WW{{\mathcal W}}
\newcommand\XX{{\mathcal X}}
\newcommand\YY{{\mathcal Y}}
\newcommand\ZZ{{\mathcal Z}}
\newcommand\CH{{\CC\HH}}
\newcommand\TC{{\TT\CC}}
\newcommand\EXH{{ \EE (X, \HH )}}
\newcommand\GXH{{ \GG (X, \HH )}}
\newcommand\GYH{{ \GG (Y, \HH )}}
\newcommand\PEX{{\PP\EE  (X, \HH , \GG , \LL )}}
\newcommand\MF{{\MM\FF}}
\newcommand\PMF{{\PP\kern-2pt\MM\FF}}
\newcommand\ML{{\MM\LL}}
\newcommand\PML{{\PP\kern-2pt\MM\LL}}
\newcommand\GL{{\GG\LL}}
\newcommand\Pol{{\mathcal P}}
\newcommand\half{{\textstyle{\frac12}}}
\newcommand\Half{{\frac12}}
\newcommand\Mod{\operatorname{Mod}}
\newcommand\Area{\operatorname{Area}}
\newcommand\ep{\epsilon}
\newcommand\hhat{\widehat}
\newcommand\Proj{{\mathbf P}}
\newcommand\U{{\mathbf U}}
 \newcommand\Hyp{{\mathbf H}}
\newcommand\D{{\mathbf D}}
\newcommand\Z{{\mathbb Z}}
\newcommand\R{{\mathbb R}}
\newcommand\Q{{\mathbb Q}}
\newcommand\E{{\mathbb E}}
\newcommand\til{\widetilde}
\newcommand\length{\operatorname{length}}
\newcommand\tr{\operatorname{tr}}
\newcommand\gesim{\succ}
\newcommand\lesim{\prec}
\newcommand\simle{\lesim}
\newcommand\simge{\gesim}
\newcommand{\simmult}{\asymp}
\newcommand{\simadd}{\mathrel{\overset{\text{\tiny $+$}}{\sim}}}
\newcommand{\ssm}{\setminus}
\newcommand{\diam}{\operatorname{diam}}
\newcommand{\pair}[1]{\langle #1\rangle}
\newcommand{\T}{{\mathbf T}}
\newcommand{\inj}{\operatorname{inj}}
\newcommand{\pleat}{\operatorname{\mathbf{pleat}}}
\newcommand{\short}{\operatorname{\mathbf{short}}}
\newcommand{\vertices}{\operatorname{vert}}
\newcommand{\collar}{\operatorname{\mathbf{collar}}}
\newcommand{\bcollar}{\operatorname{\overline{\mathbf{collar}}}}
\newcommand{\I}{{\mathbf I}}
\newcommand{\tprec}{\prec_t}
\newcommand{\fprec}{\prec_f}
\newcommand{\bprec}{\prec_b}
\newcommand{\pprec}{\prec_p}
\newcommand{\ppreceq}{\preceq_p}
\newcommand{\sprec}{\prec_s}
\newcommand{\cpreceq}{\preceq_c}
\newcommand{\cprec}{\prec_c}
\newcommand{\topprec}{\prec_{\rm top}}
\newcommand{\Topprec}{\prec_{\rm TOP}}
\newcommand{\fsub}{\mathrel{\scriptstyle\searrow}}
\newcommand{\bsub}{\mathrel{\scriptstyle\swarrow}}
\newcommand{\fsubd}{\mathrel{{\scriptstyle\searrow}\kern-1ex^d\kern0.5ex}}
\newcommand{\bsubd}{\mathrel{{\scriptstyle\swarrow}\kern-1.6ex^d\kern0.8ex}}
\newcommand{\fsubeq}{\mathrel{\raise-.7ex\hbox{$\overset{\searrow}{=}$}}}
\newcommand{\bsubeq}{\mathrel{\raise-.7ex\hbox{$\overset{\swarrow}{=}$}}}
\newcommand{\tw}{\operatorname{tw}}
\newcommand{\base}{\operatorname{base}}
\newcommand{\trans}{\operatorname{trans}}
\newcommand{\rest}{|_}
\newcommand{\bbar}{\overline}
\newcommand{\UML}{\operatorname{\UU\MM\LL}}
\newcommand{\EL}{\mathcal{EL}}
\newcommand{\tsum}{\sideset{}{'}\sum}
\newcommand{\tsh}[1]{\left\{\kern-.9ex\left\{#1\right\}\kern-.9ex\right\}}
\newcommand{\Tsh}[2]{\tsh{#2}_{#1}}
\newcommand{\qeq}{\mathrel{\approx}}
\newcommand{\Qeq}[1]{\mathrel{\approx_{#1}}}
\newcommand{\qle}{\lesssim}
\newcommand{\Qle}[1]{\mathrel{\lesssim_{#1}}}
\newcommand{\simp}{\operatorname{simp}}
\newcommand{\vsucc}{\operatorname{succ}}
\newcommand{\vpred}{\operatorname{pred}}
\newcommand\fhalf[1]{\overrightarrow {#1}}
\newcommand\bhalf[1]{\overleftarrow {#1}}
\newcommand\sleft{_{\text{left}}}
\newcommand\sright{_{\text{right}}}
\newcommand\sbtop{_{\text{top}}}
\newcommand\sbot{_{\text{bot}}}
\newcommand\sll{_{\mathbf l}}
\newcommand\srr{_{\mathbf r}}
\newcommand\geod{\operatorname{\mathbf g}}
\newcommand\mtorus[1]{\boundary U(#1)}
\newcommand\A{\mathbf A}
\newcommand\Aleft[1]{\A\sleft(#1)}
\newcommand\Aright[1]{\A\sright(#1)}
\newcommand\Atop[1]{\A\sbtop(#1)}
\newcommand\Abot[1]{\A\sbot(#1)}
\newcommand\boundvert{{\boundary_{||}}}
\newcommand\storus[1]{U(#1)}
\newcommand\Momega{\omega_M}
\newcommand\nomega{\omega_\nu}
\newcommand\twist{\operatorname{tw}}
\newcommand\modl{M_\nu}
\newcommand\MT{{\mathbb T}}
\newcommand\Teich{{\mathcal T}}
\renewcommand{\Re}{\operatorname{Re}}
\renewcommand{\Im}{\operatorname{Im}}

\title{Positive area and inaccessible fixed points for hedgehogs}

\author{Kingshook Biswas }

\date{}

\thanks{Research partly supported by  Department of Science and Technology research project
grant DyNo. 100/IFD/8347/2008-2009}

\begin{abstract} Let $f$ be a germ of holomorphic
diffeomorphism with an irrationally indifferent fixed point at the
origin in $\bbC$ (i.e. $f(0) = 0, f'(0) = e^{2\pi i \alpha},
\alpha \in \bbR - \bbQ$). Perez-Marco showed the existence of a
unique family of nontrivial invariant full continua containing the
fixed point called Siegel compacta. When $f$ is non-linearizable
(i.e. not holomorphically conjugate to the rigid rotation
$R_{\alpha}(z) = e^{2\pi i \alpha}z$) the invariant compacts
obtained are called hedgehogs. Perez-Marco developed techniques
for the construction of examples of non-linearizable germs; these
were used by the author to construct hedgehogs of Hausdorff
dimension one, and adapted by Cheritat to construct Siegel disks
with pseudo-circle boundaries. We use these techniques to
construct hedgehogs of positive area and hedgehogs with inaccessible fixed
points.

\smallskip

\begin{center}


\end{center}

\end{abstract}

\maketitle

\overfullrule=0pt

\tableofcontents

\section{Introduction}

\bigskip

A germ $f(z) = e^{2\pi i \alpha}z + O(z^2), \alpha \in \mathbb{R -
Q}/\mathbb{Z}$ of holomorphic diffeomorphism fixing the origin in
$\mathbb{C}$ is said to be linearizable if it is analytically
conjugate to the rigid rotation $R_{\alpha}(z) = e^{2\pi i \alpha}
z$. The number $\alpha$ is called the rotation number of $f$, and
the maximal domain of linearization is called the Siegel disk of
$f$. The linearizability of $f$ is dependent on the arithmetic of
$\alpha$, and the optimal arithmetic condition for linearizability
is known in this setting (see \cite{siegel}, \cite{brjuno},
\cite{yoccoz}).
%Siegel \cite{siegel} showed that if $\alpha$ satisfies a
%Diophantine condition $$ \left|\alpha - \frac{p}{q}\right| \geq
%\frac{C}{q^\tau} $$ for all integers $p,q$ for some $C, \tau >0$,
%then $f$ is linearizable. Brjuno \cite{brjuno} showed that
%linearizability holds under the more general Brjuno condition
%defined as follows: Let $(\alpha^{(i)})_{i \geq 0}$ denote the
%fractions arising from the continued fraction algorithm (given by
%iterating the Gauss map $\alpha \mapsto 1/\alpha (mod 1)$) applied
%to $\alpha$, then the Brjuno condition is $$ \mathcal{B}(\alpha)
%:= \sum_{n = 0}^{\infty} {\alpha}^{(0)} \dots {\alpha}^{(n-1)}
%\log \left(\frac{1}{{\alpha}^{(n)}}\right) < +\infty. $$
%Yoccoz \cite{yoccoz} %cite Yo
%proved the optimality of the Brjuno condition: if
%$\mathcal{B}(\alpha) = +\infty$ then there exists a
%nonlinearizable germ $f$ with rotation number $\alpha$.

\medskip

In this article we are primarily concerned with non-linearizable
germs. Perez-Marco studied the dynamics of such germs,
%(\cite{perezmens}, \cite{perezmtopology}, \cite{perezmdynamics},
%\cite{perezmcirclemaps}), \cite(perezminvent}),
proving the existence of non-trivial invariant full continua containing the
fixed point called {\it hedgehogs} (\cite{perezmcirclemaps}).
These have empty interior and are not locally connected at any
point except perhaps the fixed point (\cite{perezmtopology},
\cite{perezmdynamics}). In \cite{perezmens}, \cite{perezminvent},
\cite{perezmsmooth}, Perez-Marco developed techniques using  "tube-log Riemann
surfaces" for the construction of interesting examples of
indifferent dynamics. These were used by the author to construct
examples of hedgehogs containing smooth combs and hedgehogs of
minimal Hausdorff dimension one (\cite{biswas1}, \cite{biswas2}),
and further developed by Cheritat to construct Siegel disks with
pseudo-circle boundaries (\cite{cheritat}).

\medskip

We use these techniques, incorporating Cheritat's modifications,
to construct hedgehogs of positive area:

\medskip

\begin{theorem} \label{parea}{There exists a non-linearizable germ
with a hedgehog of positive area.}
\end{theorem}

\medskip

As in the previous constructions of hedgehogs, the
non-linearizable germ is obtained as a limit of a sequence of
finite order germs which have decreasing linearization domains
nesting down to a hedgehog of the limiting germ. Thanks to the
modified construction introduced by Cheritat, these linearization
domains can be constructed such that the loss of area at each step
is as small as desired, ensuring that the intersection has
positive area. While the above construction gives univalent maps
with hedgehogs of positive area, it is an interesting question
whether this can be achieved in the setting of entire functions,
in particular polynomials. For quadratic polynomials with a
nonlinearizable irrationally indifferent fixed point (or {\it
Cremer point}), the hedgehogs belong to the Julia set, and Buff
and Cheritat have constructed such Cremer quadratic polynomials
with positive area Julia sets (\cite{buffcheritat}); it is unclear
however whether the corresponding hedgehogs have positive area.

\medskip

%
%\medskip
%
%We can also construct similar examples for linearizable maps
%extending to a neighbourhood of the Siegel disk:
%
%\medskip
%
%\begin{theorem} \label{dimtwodisk} There exists a linearizable
%germ $f$ univalent on a neighbourhood of the closure of its Siegel
%disk $D$ such that $\partial D$ has Hausdorff dimension two.
%\end{theorem}
%
%\medskip

In \cite{perezmtopology} Perez-Marco posed the question of whether
the fixed point of a hedgehog is always accessible from its
complement. The answer is negative:

\medskip

\begin{theorem} \label{inaccessible} There exists a
nonlinearizable germ $f$ with a hedgehog $K$ such that
the fixed point is inaccessible from $\hat{\bbC} -
K$.
\end{theorem}

\medskip

The above example is constructed controlling the geometry of the complements
of the decreasing linearization domains, in such a way that the fixed point becomes inaccessible
from the union of these complementary domains.

\medskip

The general scheme of the construction is described
in Section 2, and then applied to construct examples in Section 3.

\medskip

\noindent{\bf Acknowledgements} The author thanks Ricardo
Perez-Marco and Arnaud Cheritat for helpful discussions.

\medskip

\section{Outline of the construction}

\bigskip

Let $E(z) = e^{2\pi i z}$ be the universal covering $E: \mathbb{C}
\to \mathbb{C}^*$ with deck transformation $T(z) = z+1$. We will
denote by ${\bbD}_R$ the disk $\{ |z| < R \}$ and by ${\bbH}_M$
the upper half-plane $\{ \Im z > M \}$. For $\alpha > 0$ we let
$M_{\alpha}(z) = \alpha z$. Let ${\SSS}_{\alpha}$ denote the set
of univalent maps $f$ on the unit disk $\mathbb{D}$ such that
$f(0) = 0, f'(0) = e^{2\pi i \alpha}, \alpha \in \mathbb{R}$, and
$\hat{\SSS}_{\alpha}$ be the set of lifts $F$ of maps in
$\SSS_{\alpha}$ to $\mathbb{H}$, satisfying $F(z) = z + \alpha +
\phi(z)$, where $\phi$ is 1-periodic and and tends to $0$ as $\Im
z \to +\infty$. We will refer to the $T$-invariant preimages by
$E$ of Jordan domains in $\bbC$ containing the origin as
$1$-periodic domains. By {\it linearization domain (respectively
hedgehog)} we will refer to both a linearization domain
(respectively hedgehog) $K \subset \bbD$ for $f$ as well the
$T$-invariant, $F$-invariant set $E^{-1}(K) \subset \mathbb{H}$.
By {\it linearizable maps} we will refer to both maps $f \in
\SSS_{\alpha}$ conjugate to rotations and their lifts $F \in
\hat{\SSS}_{\alpha}$. A {\it semi-flow} $\mathcal{F}$ of maps will
refer to a commuting collection of maps $(f_t)_{t \in A}$ such
that $f_t \in \SSS_t$ (or maps $(F_t)_{t \in A}, F_t \in
\hat{\SSS}_t$), for some $A \subset \bbR$; $\mathcal{F}$ will be
called linearizable if all maps in $\mathcal{F}$ are linearizable.

\medskip

\subsection{Approximation by linearization domains}

\medskip

The key to the construction is the following approximation lemma,
due to Cheritat \cite{cheritat}, allowing one to approximate any
Jordan domain by an invariant domain for a vector field such that the elements
of the flow for small times extend univalently to a given large
disk. We let $Y$ be the vector field $+1\frac{\partial}{\partial
z}$. Any holomorphic map $f$ with nonvanishing derivative on a
domain $U$ defines a nonvanishing vector field $f_* Y$ on $U$
pulling back $Y$ by $f$.

\medskip

\begin{prop} \label{approxn} Let $D \subset \bbC$ be a Jordan domain containing
the origin, $\Lambda = E^{-1}(D)$ a $1$-periodic domain and $\epsilon, M > 0$.
There exists $\delta = \delta(\Lambda, \epsilon, M) >
0$, an entire function $\Psi = \Psi(\Lambda, \epsilon, M)$ with nonvanishing
derivative of the form $\Psi(z) = z - iC + \psi(z)$ where $C > 0$,
$\psi$ is $1$-periodic, $\psi(z) \to 0$ as $\Im z \to +\infty$,
and a $1$-periodic domain $\Omega = \Omega(\Lambda, \epsilon, M) \subset \Lambda$
such that the following holds for the vector field $X = \Psi_* Y$:

\smallskip

$\Psi$ maps $\Omega$ univalently to the upper half-plane $\bbH$,
$\partial \Omega$ is contained in the
$\epsilon$-neighbourhood of $\partial \Lambda$,
and for $|t| < \delta$, the maps $F_t = exp(tX)$ are univalent on
$\bbH_{-M}$ and leave $\Omega$ invariant.
\end{prop}

\medskip

\noindent{\bf Proof:} Let $\phi : \overline{D} \to \overline{\bbD}$ be the
normalized Riemann mapping such that $\phi(0) = 0, \phi'(0) > 0$. Let
$\Phi : \overline{\Lambda} \to \overline{\bbH}$ be a lift of
$\phi$ to $\Lambda$, so that $\Phi$ commutes with $T$, and
$\Phi(z) = cz + \chi(z)$ where $\chi(z) = O(e^{2\pi i z})$ is 1-periodic.
Then $\Phi'$ is non-zero, 1-periodic and hence can be written in
the form $\Phi' = ce^{G}$ where $G$ is 1-periodic, holomorphic on
$\Lambda$ and tends to $0$ as $\Im z \to +\infty$, and is therefore of
the form $G = g \circ E$ for some $g$ holomorphic on $D$ with $g(0) = 0$.

\medskip

Fix $\kappa > 0$. By Runge's Theorem we can approximate $g$ uniformly by a polyomial
$P$ on any given compact subdomain of $D$ such that $P(0) = 0$; it follows easily that we
can approximate $\Phi$ uniformly on $\Lambda_{\kappa} = \{ \Im
\Phi(z) \geq \kappa \}$ by an entire function $\Psi_1$ commuting
with $T$ with nowhere vanishing derivative $\Psi'_1(z) = ce^{P \circ E(z)}$ on $\bbC$,
such that $\Psi_1(z) = cz + \psi(z)$ where $\psi$ is $1$-periodic and tends to $0$
as $\Im z \to +\infty$. Choosing the
approximation close enough, we may assume that $\Psi_1$ is univalent
on $\Lambda_{\kappa}$ and $\Psi_1(\Lambda_{3\kappa}) \subset \Phi(\Lambda_{2\kappa}) \subset
\Psi_1(\Lambda_{\kappa})$, so that letting $\Omega = \{ z \in \Lambda : \Im \Psi_1(z) > 2\kappa\}$,
we have $\Lambda_{3\kappa} \subset \Omega \subset
\Lambda_{\kappa}$. Let $\Psi = \Psi_1 - 2i\kappa$, then $\Psi$
maps $\Omega$ univalently onto $\bbH$.

\medskip

The pull-back $X = \Psi_* Y$ is an entire non-zero 1-periodic vector field.
Hence there exists $\delta > 0$ such that the maps
$F_t = exp(tX), |t| < \delta$, are defined and univalent on the upper half-plane $\bbH_{-M}$.
Since the flow of $Y$ leaves upper half-planes invariant,
the maps $(F_t)_{|t| < \delta}$ leave the domain $\Omega$
invariant. Taking $\kappa$ small enough
initially, we may ensure $\partial \Omega$ is within the
$\epsilon$-neighbourhood of $\partial \Lambda$. $\diamond$

\medskip

\subsection{Iterated inverse renormalization}

%\medskip
%
%We can perform an iterated inverse renormalization to obtain a sequence
%of semi-flows while controlling the geometry of the limiting invariant compacts as
%follows:

\medskip

Fix a strictly decreasing sequence constants $M_n > 1, n \geq 0$.
Let $\Lambda_n \subset \bbH, n \geq 0$ be a given sequence of
$1$-periodic domains and $\epsilon_n > 0$ a given sequence of positive
reals. We then construct inductively a sequence of corresponding semi-flows $\mathcal{F}_n$ as
follows below; the flexibility of the construction derives from the fact
that in practice the domains $\Lambda_n$ can be
chosen appropriately to obtain examples with specified properties.

\medskip

We then construct inductively the following:

\medskip

\begin{itemize}

\item A sequence $(\Psi_n)_{n \geq 0}$ of entire functions
with nonvanishing derivatives of the form $\Psi_n(z) = z - iC_n + \psi_n(z)$ where
$C_n > 0$, $\psi_n$ is $1$-periodic and tends to $0$ as $\Im z \to +\infty$.

\item A sequence of $1$-periodic domains $(\Omega_n)_{n \geq 0}$
such that $\Psi_n$ maps $\Omega_n$ univalently onto $\bbH$.

\item A sequence $\delta_n > 0, n \geq 0$.

\item A sequence of integers $a_n \geq 3, n \geq 0$.

\item Finite subsets $A_n \subset \bbQ$ and linearizable semi-flows
$\mathcal{F}_{n} = (F_{n,t})_{t \in A_n}$ for $n \geq 0$ such that
the sets $A_n$ are increasing and the semi-flows $\mathcal{F}_n$ are univalent
on $\bbH_{-M_n}$.

\item A decreasing sequence of $1$-periodic domains $D_n \subset \bbH, n \geq 0$
such that $D_n$ is invariant under $\mathcal{F}_n$.
\end{itemize}

\medskip

At stage $n = 0$: Given $\Lambda_0, \epsilon_0, M_0$ we let $\Psi_0 = \Psi(\Lambda_0,
\epsilon_0, M_0), \Omega_0 = \Omega(\Lambda_0, \epsilon_0, M_0)$ and
$\delta_0 = \delta(\Lambda_0, \epsilon_0, M_0)$ be as
given by Proposition \ref{approxn}, and $X_0 = (\Psi_0)_* Y$. Let $a_0
\geq 3$ be an integer such that $1/a_0 < \delta_0$. We let
$A_{0} = \{ \varepsilon / a_0 : \varepsilon \in \{-1,0,1\} \}$ and
$\mathcal{F}_0 = (exp(tX_0))_{t \in A_0}$. We let $D_0 = \Omega_0$. We assume $a_0$ is
chosen large enough so that the maps $exp(tX_0)$ map $\bbH_{-M_1}$
univalently into $\bbH_{-M_0}$ for $|t| \leq 1/a_0$.

\medskip

At stage $n \geq 1$: We assume by induction that we have the
following induction hypotheses:

\medskip

\begin{enumerate}

\item The maps $F_{n-1,t}$ in $\mathcal{F}_{n-1}$ map $\bbH_{-M_n}$ univalently into
$\bbH_{-M_{n-1}}$.

\item The map $\Phi_{n-1} := M_{a_{n-1}} \circ \Psi_{n-1} \circ \dots \circ M_{a_0}
\circ \Psi_0$ maps $D_{n-1}$ univalently onto $\bbH$.

\item $\Phi_{n-1}$ semi-conjugates the maps in $\mathcal{F}_{n-1}$ to powers of $T$.

\end{enumerate}

\medskip

Each $\Psi_k$ is of the form $\Psi_k = f_k \circ E$ for some entire function $f_k$, hence
the image of $\bbH_{-M_0}$ under $\Phi_n$ is contained in a half-plane
$\bbH_{-\hat{M}_n}$ for some $\hat{M}_n > 0$. We let
$\Psi_n = \Psi(\Lambda_n, \epsilon_n, \hat{M}_n),
\Omega_n = \Omega(\Lambda_n, \epsilon_n, \hat{M}_n)$ and
$\delta_n = \delta(\Lambda_n, \epsilon_n, \hat{M}_n)$  be as
given by Proposition \ref{approxn}. Let $\hat{X}_n = (\Psi_n)_* Y$
and $X_n = (\Phi_{n-1})_* \hat{X}_n$.
For $|t| < \delta_n$ the maps $exp(t\hat{X}_n)$ are univalent on $\bbH_{-\hat{M}_n}$,
hence for $|t| < \delta_n$ the maps $exp(tX_n)$ are univalent on
$\bbH_{-M_0}$. Moreover $\Phi_{n-1}$ semi-conjugates the maps
$exp(tX_{n-1}), t \in A_{n-1}$ to powers of $T$, which commute
with the maps $exp(t\hat{X}_n), |t| < \delta_n$, hence the maps
$exp(tX_n), |t| < \delta_n$, commute with the maps $exp(tX_{n-1}),
t \in A_{n-1}$.

\medskip

We let $a_n \geq 3$ be an integer such that $1/a_n < \delta_n$. We
let
\begin{align*}
A_n & = \left\{ t = t_1 + \frac{t_2}{a_0 \dots a_{n-1}} : t_1 \in A_{n-1}, t_2 \in \{-1/a_n,0,1/a_n\}
\right\}, \\
\mathcal{F}_n & = \left( F_{n,t} = F_{n-1,t_1} \circ exp(t_2 X_n): t = t_1+\frac{t_2}{a_0 \dots a_{n-1}}
\in A_n\right)_{t \in A_n}. \\
\end{align*}

We assume that $a_n$ is chosen large enough so that the maps
$exp(t X_n)$ map $\bbH_{-M_{n+1}}$ univalently into $\bbH_{-M_n}$ for $|t| \leq
1/a_n$. Then the induction hypotheses (1),(3) above are satisfied. Finally, $\Phi_{n-1}$ maps
$D_{n-1}$ univalently onto $\bbH$, and $M_{a_n} \circ \Psi_n$ maps
$\Omega_n \subset \bbH$ univalently onto $\bbH$. We let
$D_n = \Phi^{-1}_{n-1}(\Omega_n)$. Then induction hypothesis (2) is satisfied and the
induction can proceed, giving a sequence of semi-flows
$\mathcal{F}_n$.

\medskip

Let $A = \cup_{n \geq 0} A_n$. We observe that $F_{n,t} = F_{n+1,t}$ for $t \in A_n \cap
A_{n+1}$ so the semi-flows $\mathcal{F}_n$ are increasing, giving a well-defined
semi-flow $\mathcal{F} = (F_t = F_{n,t} : t \in A)_{t \in A}$ univalent on $\bbH_{-1}$
(since $M_n > 1$ for all $n$). The maps in $\mathcal{F}$ are lifts of maps in a semi-flow
$\overline{\mathcal{F}} = (f_t)_{t \in A}$ univalent in the disc
$\bbD_{e^{-2\pi}}$ such that $f_t(0) = 0, f'_t(0) = e^{2\pi i t}$.
Let $K_n = \overline{E(D_n)}$ and

\begin{align*}
A_{\infty} & = \overline{A} = \left\{ \sum_{n = 0}^{\infty}
\frac{\varepsilon}{a_0 \dots a_n} : \varepsilon \in \{-1,0,1\}
\right\} \\
D_{\infty} & = \overline{\cap_{n \geq 0} D_n} \\
K_{\infty} & = E(D_{\infty}) \\
\end{align*}

\begin{lemma} \label{limitflow} There are unique semi-flows $\mathcal{F}_{\infty} =
(F_t)_{t \in A_{\infty}}$ (respectively $\overline{\mathcal{F}}_{\infty} = (f_t)_{t \in A_{\infty}}$)
containing $\mathcal{F}$ (respectively containing
$\overline{\mathcal{F}})$. The semi-flows $\mathcal{F}_{\infty}$
(respectively $\overline{\mathcal{F}}_{\infty}$) leave invariant
$D_{\infty}$ (respectively $K_{\infty}$).
\end{lemma}

\medskip

\noindent{\bf Proof:} Let $t_n \in A, t_n \to t \in A_{\infty}$. The maps $F_{t_n}$ are
lifts of maps $f_{t_n}$ univalent in $\bbD$ such that $f_{t_n}(0) = 0, f'_{t_n}(0) = e^{2\pi i
t_n}$ and hence form a normal family. Any subsequential limit $F$ of $F_{t_n}$ is a lift
of a map $f$ univalent in $\bbD$ with $f(0) = 0, f'(0) = e^{2\pi i t}$; moreover $f$ commutes with
all the maps in $f_{t_n}$. It is easy to see
that these two properties determine the Taylor expansion of $f$ at the origin
(and hence $f$) uniquely (because $A$
accumulates at $0$). Consequently $F$ is the unique lift of $f$ tangent to the translation
$z \mapsto z+t$ as $\Im z \to +\infty$. Hence for any $t_n \to t$ the sequence $F_{t_n}$ has a
unique limit depending only on $t$, say $F_t$. We thus obtain
semi-flows $\mathcal{F}_{\infty} := (F_t)_{t \in A_{\infty}},
\overline{\mathcal{F}}_{\infty} := (f_t)_{t \in A_{\infty}}$ as
required. Since $D_m$ is invariant under $\mathcal{F}_n$
for $m \geq n$, it follows that $D_{\infty}$ is
invariant under $\mathcal{F}_{\infty}$. $\diamond$

\medskip

We note that in the above construction at each stage $n$ the
choices of $\Lambda_n, \epsilon_n$ are independent of the objects constructed up
to stage $n-1$, hence the approximating domains $\Omega_n$ can be
constructed almost arbitrarily. This allows us to control the
geometry of the limiting invariant compact $K_{\infty}$. If $K_{\infty} \neq \{0\}$,
then $\overline{\mathcal{F}}_{\infty}$ will be linearizable if and
only if the origin is an interior point of $K_{\infty}$ ; otherwise all
irrational elements of $\overline{\mathcal{F}}_{\infty}$ will be linearizable
and will have $K_{\infty}$ as a common hedgehog.
%
%Let $\Psi_0$ be an entire function such that $\Psi'_0$ is nonvanishing,
%$1$-periodic, and $\Psi_0$ maps a Jordan domain $\Omega_0 \subset
%\bbH_{-M_0}$ univalently to the upper-half plane $\bbH$. Let
%$\hat{X}_0$ be the pullback of $\hat{Y}$ by $\Psi_0$ and
%$\delta_0 > 0$ be such that the maps $exp(t\hat{X}_0)$ are
%univalent on $\bbH_{-M_0}$. Let $\mathcal{F}_0 =
%(exp(t\hat{X}_0))_{|t| < \delta_0}$.
%
%\medskip
%
%We construct inductively a sequence of linearizable semi-flows $\mathcal{F}_n$
%univalent on $\bbH_{-M_n}$ as follows:
%
%\medskip
%
%Let $n \geq 0$. We suppose we are given entire functions $\Psi_k$  and
%Jordan domains $\Omega_k \subset \bbH$ for $k = 0,\dots,n$, such that $\Psi'_k$ is nonvanishing,
%$1$-periodic, and $\Psi_k$ maps $\Omega_k$ univalently to $\bbH$.
%Let $\hat{X}_n$ be the pull-back under $\Psi_0 \circ \dots \circ \Psi_{n}$
%of $\hat{Y}$ and $\delta_n > 0$ be such that for $|t| < \delta_n$
%the maps $exp(t\hat{X}_n)$ are univalent on $\bbH_{-M_n}$, leaving
%invariant $(\Psi_0 \circ \dots \circ \Psi_n)^{-1}(\bbH)$.
%
%\medskip
%
%Let $a_{n+1} \geq 2$ be an integer larger than $1/\delta_n$. By
%Proposition \ref{approxn} we obtain an entire function
%$\Psi_{n+1}$ with $1$-periodic nonvanishing derivative and a
%domain $\Omega_{n+1} \subset \Lambda_{n+1}$ such that $\Psi_{n+1}$ maps
%$\Omega_{n+1}$ univalently to $\bbH$ and $\partial \Omega_{n+1}$
%is within the $\epsilon_{n+1}$ neighbourhood of $\partial
%\Lambda_{n+1}$.
%

\bigskip

\section{Construction of examples}

\bigskip

Varying the choices of the domains $\Lambda_n$ in the above
construction will allow us to construct the required examples. We
fix a sequence $M_n > 1$ as in the previous section.

\bigskip

\subsection{Positive area hedgehogs}

\bigskip

Fix $0 < r < 1$ and sequences $0 < r_n < r, \alpha_n < 1$ such
that $r_n \to 0$ and $\prod_{n \geq 0} \alpha_n
> 0$.

\medskip

\begin{lemma} \label{areachoice} We can choose the
domains $\Lambda_n$ and constants $\epsilon_n > 0$
inductively in order that the Jordan domains $K_n = \overline{E(D_n)}$
given by the construction satisfy the following:

\begin{enumerate}

\item $\partial K_n$ intersects $\{ |z| > r \}$ for all $n \geq 0$.

\item $\partial K_n$ intersects $\cap \{|z| < r_n \}$ for all $n \geq 0$.

\item $\lambda(K_n)/\lambda(K_{n-1}) > \alpha_n$ for all $n \geq 1$ where $\lambda$ is
Lebesgue measure on $\bbC$.

\end{enumerate}
\end{lemma}

\medskip

\noindent{\bf Proof:} In the construction of the previous section,
the domains $K_n$ are given by
$K_n = \overline{E(D_n)}$ where $D_n = \Phi^{-1}_{n-1}(\Omega_n)$.
The domains $\Omega_n$ are given by Proposition \ref{approxn} as
approximations to the domains $\Lambda_n$.
At each stage $n$ we will choose
$\Lambda_n$ such that the compact $K'_n =
\overline{E(D'_n)}, D'_n = \Phi^{-1}_{n-1}(\Lambda_n)$, satisfies
conditions (1),(2) above (with $K_n$ replaced by $K'_n$), and condition (3)
with $K_n$ replaced by $K'_n$, namely $\lambda(K'_n) / \lambda(K_{n-1}) > \alpha_n$.
This will suffice as the approximating domain $\Omega_n$ can then be chosen
close enough so that the compact $K_n$ satisfies (1), (2), (3).

\medskip

Let $h_n > h > 0$ be such that $e^{-2\pi h} = r, e^{-2\pi h_n} = r_n$.
For $H > 0$ let $S(H)$ be the open strip $\{ 0 < \Im z < H\}$.

\medskip

At stage $n = 0$: We choose $\Lambda_0 \subset \bbH$ a $1$-periodic domain
such that $\partial \Lambda_0$ intersects $S(h)$ and $\bbH_{h_0}$.
Then conditions (1), (2) are satisfied by $K'_0 =
E(\Lambda_0)$ (condition (3) is empty for $n = 0$).
%By Proposition \ref{approxn}, we can choose $\epsilon_0 >
%0$ small enough so that $\partial \Omega_0 \cap l(h) \neq
%\emptyset, \partial \Omega_0 \cap l(h_n) \neq \emptyset$.

\medskip

%We assume that $\Lambda_k, \epsilon_k, 0 \leq k \leq
%n-1$, have been chosen, and hence $\Psi_k, a_k, 0 \leq k \leq n-1$
%have been determined, such that conditions (1),(2),(3) above hold for $K_k$
%and $(1),(2)$ hold for $K'_k$ for  $0 \leq k \leq n-1$.
%
%\medskip

At stage $n \geq 1$: If we choose $\Lambda_n$ such that $\partial \Lambda_n$ intersects the
images of the nonempty open sets $S(h) \cap D_{n-1}$ (nonempty by induction) and $\bbH_{h_n} \cap D_{n-1}$
under $\Phi_{n-1}$, then conditions (1), (2) will be satisfied by
$K'_n$. Since $K_{n-1} = \overline{E(\Phi^{-1}_{n-1}(\bbH))}$, if we
choose $\Lambda_n$ to "fill out" most of $\bbH$ we can ensure that the Jordan domain
$K'_n = E(\Phi^{-1}_{n-1}(\Lambda'_n))$ satisfies condition (3), $\lambda(K'_n) /
\lambda(K_{n-1}) > \alpha_n$. It is clear that $\Lambda_n$
can be chosen to satisfy these constraints. This completes the
induction. $\diamond$

\medskip

\noindent{\bf Proof of Theorem \ref{parea} :} Since $\partial K_n$
intersects $\{z = r\}$ and $\{|z| = r_n\}$ for all $n$ it follows that
$K_{\infty} \neq \{0\}$ and the origin is not an
interior point of $K_{\infty}$, therefore all
elements $f_t$ of the semi-flow $\overline{\mathcal{F}}_{\infty}$ are
nonlinearizable for $t$ irrational, and $K_{\infty}$ is a common
hedgehog for these elements.

\medskip

Since the compacts $K_n$ decrease to $K_{\infty}$, $\lambda(K_n)
\to \lambda(K_{\infty})$ as $n \to \infty$, and
$$
\frac{\lambda(K_N)}{\lambda(K_0)} = \prod_{n = 1}^{N}
\frac{\lambda(K_{n})}{\lambda(K_{n-1})} \geq \prod_{n = 1}^{N} \alpha_n \to \prod_{n = 1}^{\infty} \alpha_n > 0
$$
so $\lambda(K_{\infty}) > 0$. $\diamond$

\bigskip

\bigskip

\subsection{Hedgehogs with inaccessible fixed points}

\bigskip

We give first a brief outline of the construction. We wish to choose the
domains $\Lambda_n$ inductively controlling the geometry of the
complementary $1$-periodic domains $W_n := \mathbb{H} -
\overline{D_n}$ so that there is no continuous curve $\gamma
: [0,+\infty) \to W_{\infty}$ with $\lim_{t\to +\infty} \Im
\gamma(t) = +\infty$, but $\sup_{z \in W_{\infty}} \Im z = +\infty$,
where $W_{\infty} = \cup_{n \geq 1} W_n$ is
the pre-image under $E$ of the complement $\mathbb{C} -
K_{\infty}$ of the invariant compact obtained from the
construction. Then $K_{\infty}$ will be a hedgehog such that the
origin is inaccessible from $\mathbb{C} - K_{\infty}$.

\medskip

Given $W_n$, we construct inductively $W_{n+1}$ by choosing an
equipotential curve $\gamma'_n$ of $D_n$ very close to $\gamma_n := \partial D_n =
\partial W_n$, then constructing a Jordan arc $\alpha_n$ lying
inside the "upper" subdomain of $\mathbb{C}$ bounded by $\gamma'_n$, such
that $\alpha_n$ starts from a point of $\gamma'_n$, follows
$\gamma'_n$ closely for a full period, then goes up to a height
much larger than the largest height of any point of $\gamma_n$,
and such that the images of $\alpha_n$ under the semi-flow
$\mathcal{F}_n$ are pairwise disjoint. We then choose $D_{n+1}$ such that $W_{n+1}$
is very close to the domain given by the union of the "lower"
subdomain of $\mathbb{C}$ bounded by $\gamma'_n$ together with the images
under $\mathcal{F}_n$ of a small neighbourhood of $\alpha_n$.

\medskip

We refer to these small neighbourhoods of the Jordan arcs $F(\alpha_n), F \in \mathcal{F}_n$,
as "tubes". Roughly speaking, each domain $W_{n+1}$ is obtained from the preceding domain
$W_n$ by adding a collection of tubes. Thus each $W_n$ may be
written as a union of tubes, where the tubes are divided into $n$
groups, the $j$th group containing the tubes added at level $j$ of
the construction. The tubes are constructed such that the tubes of
level $n+1$ go up to a height $h_{n+1}$ much larger than the
height $h_n$ reached by the level $n$ tubes, but the level $n+1$
tubes follow initially the tubes of the previous levels so closely that they must pass
through a fixed height $0 < C < 5$ before reaching a height
$h_{n+1}$. This property of the tubes forces any curve $\gamma$ in
$W_{\infty}$ with $\Im \gamma$ accumulating $+\infty$ to also have
$\Im \gamma$ accumulating $(-\infty, 2]$, hence there is no curve
with $\lim_{t \to +\infty} \Im \gamma(t) = +\infty$.

\medskip

This ends the outline, we give now the details of the actual
construction:

\medskip

The building block for the construction will be a Jordan arc
$\alpha = \alpha(q, \tau, H)$ depending on three parameters $q, \tau, H$,
where $q \geq 1$ is an integer,
and $\tau, H > 0$. The curve $\alpha$ is given by the straight
line segment from $0$ to $q+i\tau$, followed by the vertical segment from $q+i\tau$ to
$q+iH$. We note that the images
of $\alpha$ under the translation by one $T : z \mapsto z+1$ are pairwise disjoint. In practice, $q$
will be given, and $\tau$ and $H$ will have to be chosen small
and large enough respectively.

\medskip

We now define inductively the domains $\Lambda_n$, constants
$\epsilon_n, h_n, \tau_n, H_n > 0$, and also the tubes referred to above, such that
at each stage $n$ we have a decomposition of $W_n$ into a union of
tubes of levels $1$ to $n$, as follows:

\medskip

At stage $n = 0$: We let $H_0 = 10, \tau_0 = 1/2, \alpha_0 =
\alpha(2, \tau_0, H_0)$. We let $U_0$ be a $1/8$-neighbourhood of
$\alpha_0$ and let $\Lambda_0$ be the union of the lower
half-plane with the open sets $T^k(U_0), k \in \mathbb{Z}$. We
choose $\epsilon_0 > 0$ small so that we may express the complement $W_0$ of $D_0$ as
$$
W_0 = G_0 \cup \left( \bigsqcup_k T^k(V_0) \right)
$$
where $G_0$ is a the "lower" subdomain of $\mathbb{C}$ bounded by
a $1$-periodic curve very close to the real axis, and $V_0$ is a
Jordan domain, with boundary given by a Jordan arc following
closely the boundary of $U_0$ together with a cross-cut $C_0$ of
$D_0$ of diameter less than $2$
joining two points $z$ and $z+1$ with $z$ close to $-1/2$, as shown in the figure below:

{\hfill {\centerline {\psfig {figure=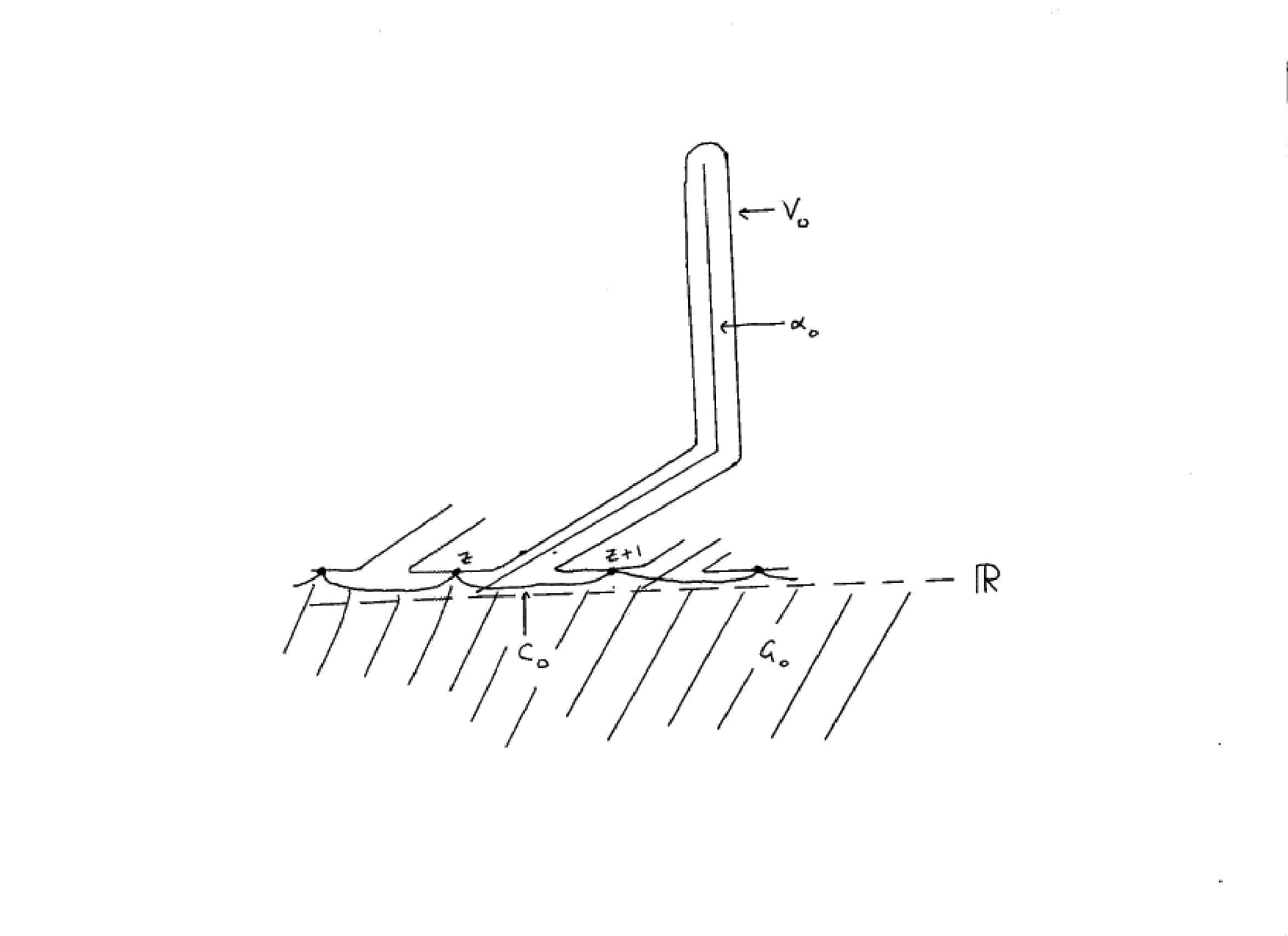,height=10cm}}}}

{\centerline {\bf Stage $n=0$ of the construction }}

\medskip

Since the choice of $\Psi_0$ is fixed by the above and
$\Phi_0 = M_{a_0} \circ \Psi_0$, we can ensure that $a_0 \geq 3$ is
chosen large enough so that $|\Phi^{-1}_0(x) - \Phi^{-1}(x')| \leq
1/2^5$ for all $x,x' \in \mathbb{R}$ with $|x - x'| \leq 2$.

\medskip

We define the height $h_0$ by $h_0 := \sup_{z \in W_0} \Im z$.

\medskip

We refer to the sets $V_{k,0,0} := T^k(V_0)$ as the level $0$
tubes of stage $0$ of the construction, and denote the bounding
cross-cuts $T^k(C_0)$ by $c_{k,0,0}$.

\medskip

At stage $n \geq 1$: We assume by induction that the $1$-periodic
domain $D_{n-1}$ (and hence also its complement $W_{n-1}$) has
been constructed, such that we have a decomposition $$ W_{n-1} =
G_{n-1} \cup \left( \bigsqcup_{j = 0}^{n-1} \bigsqcup_{k \in
\mathbb{Z}} V_{k,j,n-1} \right) $$ where $G_n$ is the "lower"
subdomain of $\mathbb{C}$ bounded by a $1$-periodic curve close to
the real axis, and $V_{k,j,n-1}$ are the level $j$ tubes of the
stage $n-1$ of the construction, with bounding cross-cuts
$c_{k,j,n-1}$.

\medskip

We assume by induction that the tubes and cross-cuts are arranged
as follows:

\medskip

Given a domain $W \subset \mathbb{C}$ bounded by a Jordan curve
$\partial W$ in $\hat{\mathbb{C}}$ with $\infty \in
\partial W$, and a cross-cut $c \subset \mathbb{C}$ of $W$, we
denote by $B(W, c)$ the bounded component of $W - c$. We assume by
induction the nesting property that for $0 \leq j \leq j' \leq n-1, k,k' \in
\mathbb{Z}$, either $B(W_{n-1}, c_{k',j',n-1}) \cap B(W_{n-1},
c_{k,j,n-1} = \emptyset$ or $B(W_{n-1}, c_{k',j',n-1}) \subset
B(W_{n-1}, c_{k,j,n-1})$. Moreover each level $j$ tube may be
written as
$$
V_{k,j,n-1} = B(W_{n-1}, c_{k,j,n-1}) - \left(\bigcup_{j' > j, k' \in
\mathbb{Z}} B(W_{n-1}, c_{k',j',n-1})\right)
$$
Each level $0$ tube $V_{k,0,n-1}$ intersects $G_{n-1}$ precisely along the cross-cut
$c_{k,0,n-1}$. Also, each tube $V_{k,j,n-1}$ of level $j \geq 1$ intersects exactly one tube
of level $j-1$ along the cross-cut $c_{k,j,n-1}$, and only
intersects tubes of levels $j' > j$ for $j' = j + 1$ along
cross-cuts of level $j+1$. We assume by induction that the heights of the tubes satisfy
$$
\sup_{z \in V_{k,j,n-1}} \Im z \leq h_j + 1/2^j + \dots 1/2^{n-1}
$$

We assume by induction that $|\Phi^{-1}_{n-1}(x) -
\Phi^{-1}_{n-1}(x')| \leq 1/2^{n+4}$ for $x,x' \in \mathbb{R}, |x
- x'| \leq 2$, and that the
interiors of the cross-cuts $c_{j,k,n-1}$ are pairwise disjoint,
and all of them are of diameter less than $1/2^{n-2}$.

\medskip

By induction we also have a univalent map $\Phi_{n-1}$ mapping
$D_{n-1}$ to $\mathbb{H}$ and conjugating $\mathcal{F}_{n-1}$ to
powers of $T$.

\medskip

We choose $\tau_n > 0$ small enough so that

$$
|\Phi^{-1}_{n-1}(x+iy) - \Phi^{-1}_{n-1}(x)| < 1/2^{n+4}, \ \ x \in \mathbb{R}, y \in [0, 10\tau_n].
$$

\medskip

We choose $H_n > 0$ large
enough so that
$$
\inf_{x \in \mathbb{R}} \Im
\Phi^{-1}_{n-1}(x+iH_n) > 2h_{n-1} + 10.
$$

We define $\alpha_{n,n}$ to be the curve
$\alpha(3 a_0 \dots a_{n-1}+2, \tau_n, H_n)$ translated by $i\tau_n$, and
choose $U_{n,n}$ to be a small neighbourhood of $\alpha_{n,n}$
such that $U_{n,n} \cap \{ 0 \leq \Re z \leq a_0 \dots a_{n-1} + 1 \}
\subset \{ 0 < \Im z < 3\tau_n \}$, and such that the translates
$U_{k,n,n} := T^k(U_{n,n})$ are pairwise disjoint. We let
$\alpha_{k,n,n} := T^k(\alpha_{n,n})$.

\medskip

We define $\Lambda_n$ to be the domain
$$
\Lambda_n := \mathbb{H} - \left( \{ \Im z < \tau_n \} \cup \left(
\bigsqcup_{k \in \mathbb{Z}} U_{k,n,n} \right) \right)
$$

The end-points $e^-_{k,j,n-1}, e^+_{k,j,n-1}$ of the cross-cuts $c_{k,j,n-1}, j=0,\dots,n-1$ of
$W_{n-1}$ lie on $\partial W_{n-1} = \partial D_{n-1}$, and are
mapped by $\Phi_{n-1}$ to points $x^-_{k,j,n-1} < x^+_{k,j,n-1}$ lying on $\mathbb{R}$, with the
property $x^+_{k,j,n-1} - x^-_{k,j,n-1} \geq a_{n-1} - 1 \geq 2$ for all $k,j$, and which can
be joined to points on $\{ \Im z = \tau_n \} \subset \partial
\Lambda_n$ by small disjoint vertical segments. Moreover for $j' \geq j, k,k' \in \mathbb{Z}$, we have the
nesting property that either
$(x^-_{k',j',n-1}, x^+_{k',j',n-1}) \cap (x^-_{k,j,n-1}, x^+_{k,j,n-1}) = \emptyset$ or
$(x^-_{k',j',n-1}, x^+_{k',j',n-1}) \subset (x^-_{k,j,n-1}, x^+_{k,j,n-1})$. We choose $\epsilon_n$
small enough so that the boundary of the domain $\Omega_n$ is close enough to
$\partial \Lambda_n$ so that the points $x^{\pm}_{k,j,n-1}$ can be
joined to points of $\partial \Omega_n$ by small almost vertical
segments $l^{\pm}_{k,j,n-1}$ lying in $\{ 0 \leq \Im z \leq 3\tau_n$. We
choose the segments such that $l^{\pm}_{k,j,n-1}$ and $l^{\pm}_{k',j',n-1}$
are disjoint if $x^{\pm}_{k,j,n-1} \neq ^{\pm}_{k',j',n-1}$, while
$l^{\pm}_{k,j,n-1} = l^{\pm}_{k',j',n-1}$ if $x^{\pm}_{k,j,n-1} =
x^{\pm}_{k',j',n-1}$.

\medskip

We also assume $\epsilon_n$ is
small enough so that we may make
almost horizontal cross-cuts $h_{k,n,n}$ of $\mathbb{H} - \Omega_n$
with disjoint interiors
lying in $\{ 0 \leq \Im z \leq 2\tau_n\}$,
each joining points of $\partial \Omega_n$ close to the points $(k-1/2)+i\tau_n$ and
$(k+1/2)+i\tau_n$, such that the Jordan domains $U'_{k,n,n} := B(\mathbb{H} - \Omega_n, h_{k,n,n})$ are close to
the neighbourhoods $U_{k,n,n}$ of the curves $\alpha_{k,n,n}$. We may assume $\epsilon_n$
is small enough so that $U'_{k,n,n} \cap \{0 \leq \Re z \leq a_0\dots a_{n-1} \} \subset
\{ 0 \leq \Im z \leq 4\tau_n \}$. The cross-cuts $h_{k,n,n}$ are
chosen also such that the starting and ending points of $h_{k,n,n}$ are the ending and
starting points respectively of $h_{k-1,n,n}$ and $h_{k+1,n,n}$,
so that the complement of the cross-cuts $h_{k,n,n}$ in
$\mathbb{H} - \Omega_n$ is the disjoint union of the domains $B(\mathbb{H} - \Omega_n,
h_{k,n,n})$ and a simply connected subdomain of $\mathbb{H} -
\Omega_n$ whose closure intersects $\partial \Omega_n$ only at the
endpoints of the cross-cuts $h_{k,n,n}$.

\medskip

Define $D_n :=
\Phi^{-1}_{n-1}(\Omega_n)$, $W_n := \mathbb{H} - D_n$. We define the tubes $V_{k,j,n}$ of level $j$ and
stage $n$ and the cross-cuts $c_{k,j,n}$, for $0 \leq j \leq n$, as follows:

\medskip

Define the cross-cuts $c_{k,n,n} := \Phi^{-1}_{n-1}(h_{k,n,n}), k \in
\mathbb{Z}$. The induction hypothesis $|\Phi^{-1}_{n-1}(x) -
\Phi^{-1}_{n-1}(x')| \leq 1/2^{n+4}$ for $x,x' \in \mathbb{R}, |x
- x'| \leq 2$ implies that the cross-cuts $c_{k,n,n}$ all have
diameter less than $1/2^{n-1}$.

\medskip

For $0 \leq j \leq n-1, k \in \mathbb{Z}$, we define first
cross-cuts $c'_{k,j,n}$ to be the union of $c_{k,j,n-1}$ with $\Phi^{-1}_{n-1}(l^-_{k,j,n-1}),
\Phi^{-1}_{n-1}(l^-_{k,j,n-1})$. We note that cross-cuts
$c'_{k,j,n}, c'_{k',j',n}$ are disjoint if $c_{k,j,n-1},
c_{k',j',n-1}$ are disjoint, while if $c_{k,j,n-1}$ and
$c_{k',j',n-1}$ intersect at an endpoint, then $c'_{k,j,n},
c'_{k',j',n}$ intersect along a small segment
$\Phi^{-1}_{n-1}(l^{\pm}_{k,j,n})$. Clearly the cross-cuts
$c'_{k,j,n}$ may thus be slightly deformed near their endpoints to
obtain cross-cuts with disjoint interiors but with the same endpoints. We define the
cross-cuts $c_{k,j,n}$ to be the cross-cuts thus obtained from the
cross-cuts $c'_{k,j,n}$. The induction hypothesis on the diameters
of the cross-cuts $c_{k,j,n-1}$ being less than $1/2^{n-2}$ implies that
the diameters of the cross-cuts $c_{j,k,n}$ are less than
$1/2^{n-1}$.

\medskip

For $0 \leq j \leq n, k \in
\mathbb{Z}$, define the tube $V_{k,j,n}$ by
$$
V_{k,j,n} := B(W_n, c_{k,j,n}) - \left(\bigcup_{j' > j, k' \in
\mathbb{Z}} B(W_n, c_{k',j',n})\right)
$$
and define $G_n := W_n - (\cup_{k,j} V_{k,j,n})$. We note that for $0 \leq j \leq n-1$,
the tubes $V_{k,j,n}$ are given by the union of $V_{k,j,n-1}$
together with a domain contained in $\Phi^{-1}_{n-1}(\{0 \leq \Im
z \leq 4\tau_n \})$, hence by the hypothesis on $\tau_n$
and the induction hypothesis on the heights
of the tubes $V_{k,j,n-1}$, the heights of the tubes $V_{k,j,n}$
satisfy
$$
\sup_{z \in V_{k,j,n}} \Im z \leq h_j + 1/2^j + \dots + 1/2^n, \ \
k \in \mathbb{Z}, 0 \leq j \leq n-1.
$$

We note that any tube $V_{k,j,n}$ for $0 \leq j \leq n-1$ is
contained in the $1$ neighbourhood of $V_{k,j,j}$, and that the
tubes $V_{k,n,n}$ are given by the pre-images of the domains
$U'_{k,n,n}$ under $\Phi^{-1}_{n-1}$.

\medskip

We define the height $h_n := \sup_{z \in W_n} \Im z$ and
note that the hypothesis on $H_n$ implies that $h_n \geq 2 h_{n-1} + 2$.

\medskip

Finally, we assume that $a_n$ is chosen large enough so that $|\Phi^{-1}_{n}(x) -
\Phi^{-1}_{n}(x')| \leq 1/2^{n+5}$ for $x,x' \in \mathbb{R}, |x
- x'| \leq 2$.

\medskip

This ends the
inductive construction.

\bigskip

We let $W_{\infty} := \cup_{n \geq 0} W_n$, observing that
$W_{\infty} = \mathbb{H} - E^{-1}(K_{\infty})$. We remark that it
is easy to show by induction, using the hypothesis on the heights
$\tau_n$, that every tube of $W_n$ intersects $\{ \Im z \leq 2 \}$.

\medskip

We note that two tubes $V_{k,j,n}, V_{k',j',n}$ can intersect
along a cross-cut $c_{k'',j'',n}$ only if $|j - j'| = 1$. For each
$n \geq 0$, let $\Gamma_n$ be the graph with vertices the tubes
$V_{k,j,n}, k \in \mathbb{Z}, 0 \leq j \leq n$ and the domain $G_n$, where two vertices
are joined by an edge if and only if the corresponding domains
intersect along a cross-cut. We define the level of a vertex to be $j$ if it
corresponds to a tube $V_{k,j,n}$ and to be $-1$ if it corresponds to $G_n$. It is clear from the above
construction that $\Gamma_n$ is a rooted tree with root the vertex
corresponding to $G_n$, and for $j \geq 0$ each vertex of level $j$ is connected
to exactly one vertex of level $j-1$, and is only
connected otherwise to vertices of level $j+1$.

\medskip

We note that any curve $\gamma : [a,b] \to W_n$ with only finitely
many transverse intersections with the union of the cross-cuts
$c_{k,j,n_1}, k \in \mathbb{Z}, 0 \leq j \leq n$, gives rise to a
finite path $\tau(\gamma)$ in the graph $\Gamma_n$, obtained from $\gamma$ by
considering the tubes of $W_n$ through which $\gamma$ passes. We
call $\tau(\gamma)$ the trace of $\gamma$ in $\Gamma_n$.

\medskip

%\begin{lemma} \label{separate} Given $n \geq n_0 \geq 1$,
%let $V_{-1}, V_0, \dots, V_m$ be a sequence of tubes of $W_n$ with strictly increasing
%levels such that $V_0$ is of level $n_0$, and let $c_i$ be the cross-cut of $W_n$ which
%is the intersection of $V_{i-1}$ and $V_i$, for $i = 0, \dots, m$.
%Then there exist cross-cuts $d_0, \dots, d_m$
%of the Jordan domains $V_0, \dots, V_m$ respectively such that $\sup_{z \in d_i} \Im z \leq 4$,
%and the subdomains $U_i$ of $V_i$ bounded by
%$c_i,d_i$ are contained in $\{ \Im z \leq h_i + 1/2^i \}$, where
%$h_i = \sup_{z \in c_i} \Im z$.
%\end{lemma}
%
%\medskip

%\noindent{\bf Proof:}

\medskip

\begin{lemma} \label{height} Let $\gamma : [0, T] \to
W_n$ be a curve such that $\Im \gamma(t) > 10$ for all $0 \leq t
\leq T$.
Then there exists $M > 10$ depending only on $\gamma(0)$
such that $\Im \gamma(T) \leq M$.
\end{lemma}

\medskip

\noindent{\bf Proof:} Let $1 \leq n_0 \leq n$ be the level of the tube of
$W_n$ containing $\gamma(0)$. We will show the statement of the lemma holds for
$M = h_{n_0} + 10$.

\medskip

We first perturb $\gamma$ slightly to a curve $\gamma_1$ such that $\gamma_1$
has only
finitely many transverse intersections with the union of the
cross-cuts $c_{k,j,n}, k \in \mathbb{Z}, 0 \leq j \leq n$, $\gamma_1(T)$ belongs to the
interior of a tube of $W_n$, $|\gamma_1(T) - \gamma(T)| < 1/2$, and such that
$\Im \gamma_1(t) > 10$ for all $0 \leq t \leq T$.
Consider the trace $\tau = \tau(\gamma_1)$ of $\gamma_1$ in
$\Gamma_{n}$.

\medskip

Then,
since $\Gamma_n$ is a tree, $\tau$ contains the unique simple
path $\alpha$ between its endpoints. Moreover $\gamma_1$ consists of arcs
lying in tubes corresponding to vertices of $\alpha$, together with
arcs lying in tubes corresponding to vertices of $\tau - \alpha$. Each latter such arc
$\beta$ starts from and returns to a bounding cross-cut $c$ of a tube $V$ corresponding to a vertex of
$\alpha$. We replace each arc $\beta$ of $\gamma_1$
by a small arc $\beta'$ with the same endpoints which is contained in the
cross-cut $c$, to obtain a curve $\gamma_2 : [0,T] \to W_n$
contained entirely in the tubes corresponding to the vertices of
$\alpha$. We then perturb $\gamma_2$ slightly near each $\beta'$, replacing a subarc of $\gamma_2$
containing $\beta'$ with endpoints in the interior of $V$ by an arc lying entirely in the interior
of $V$ with the same endpoints to obtain a curve $\gamma_3 : [0,T] \to W_n$ such that the
trace of $\gamma_3$ in $W_n$ is $\alpha$, and such that $\gamma_3(T) =
\gamma_1(T)$. Since the diameters of all cross-cuts are less than
$2$, the hypothesis on $\gamma$ implies that we may assume the choices of
$\gamma_1, \gamma_2, \gamma_3$ above have been made so that $\Im \gamma_3(t)
> 8$ for $0 \leq t \leq T$.

\medskip

Let the vertices of $\alpha$ correspond to tubes $V_1, \dots,
V_m$, and let $0 = t_0 < \dots < t_m = T$ be a partition of $[0,T]$ such that
$\gamma_3([t_{i-1}, t_i]) \subset \overline{V_i}$. We note that the level of $V_1$ is $n_0$.
Since $\alpha$ is a simple path, either the levels of the tubes
$V_i$ are strictly increasing or are strictly decreasing, or $m =
1$. If either of the latter two cases hold, then the level of
$V_m$ is less than or equal to $n_0$, hence by the hypothesis on
the height of the tubes $\Im \gamma_3(T) \leq h_{n_0}+4$ and we
are done.

\medskip

We may therefore assume that $m \geq 2$ and the levels of the tubes $V_i$ are strictly
increasing. Let $V_0$ be the unique tube of level $n_0 - 1$ connected to $V_1$. Let
$c_i$ be the cross-cut of $W_n$ which is the intersection of
$V_{i-1}$ with $V_i$, for $i = 1, \dots, m-1$.

\medskip

There is an initial subarc of the boundary of $V_1$ which follows
closely the boundary of $D_{n_0-1}$ for a full period. Since $\partial
D_{n_0 - 1}$ intersects $\{ \Im z \leq 2 \}$, it follows that we
may make a cross-cut $d_1$ of $V_1$ such that $\sup_{z \in d_1}
\Im z \leq 2+1/2^{n_0}$. Of the two Jordan
subdomains of $V_1$ bounded by $d_1$, we let $Z_1$ be the subdomain $Z_1$ whose
boundary contains $c_1$, and note that $Z_1$ is contained in $\{ \Im z \leq
h_{n_0} \}$. Since $\Im \gamma_3(t) > 8$ for all $t$, $\gamma_3$
cannot intersect $d_1$ and it follows that $\gamma_3(t_1) \in
\overline{Z_1}$, hence $c_2 \subset \partial Z_1$.

\medskip

There is an initial subarc of the boundary of $V_2$ which follows
closely the boundary of $Z_1$, hence we may make a cross-cut $d_2$
of $V_2$ close to $d_1$, such that $\sup_{z \in d_2}
\Im z \leq 2+1/2^{n_0}+1/2^{n_0+1}$, and such that of the two Jordan
subdomains of $V_2$ bounded by $d_2$, the subdomain $Z_2$ whose boundary
contains $c_2$ satisfies $Z_2 \subset \{ \Im z \leq h_{n_0} +
1/2^{n_0+1} \}$. As before, $\gamma_3$ cannot intersect $d_2$,
hence $\gamma_3(t_2) \in \overline{Z_2}$, so $c_3 \subset \partial
Z_2$.

\medskip

It follows, using the hypothesis $\Im \gamma_3(t) > 8 \ \forall t$
and arguing by induction, that we may make cross-cuts $d_1,d_2,
\dots, d_m$ of the domains $V_1,V_2,\dots,V_m$ respectively,
bounding subdomains $Z_1,Z_2,\dots,Z_m$, with $c_i \subset \partial Z_i$,
such that $\sup_{z \in d_m}
\Im z \leq 2+1/2^{n_0}+1/2^{n_0+1}+\dots +1/2^{n_0+m-1}$, and such
that $Z_m \subset \{ \Im z \leq
h_{n_0}+1/2^{n_0+1}+\dots+1/2^{n_0+m-1} \}$. Since
$\gamma_3(t_{m-1}) \in c_m \subset \partial Z_m$, and $\gamma_3$
cannot intersect $d_m$, we must have $\gamma_3(T) = \gamma_3(t_m)
\in \overline{Z_m}$.

\medskip

It follows that $\Im \gamma(T) \leq \Im \gamma_3(T) + 1 \leq h_{n_0} + 1+1 < M$.
$\diamond$

\medskip

\noindent{\bf Proof of Theorem \ref{inaccessible}}: Since $D_n \cap \{ \Im z \leq 2 \} \neq \emptyset$
for all $n$ the compact $K_{\infty}$ is non-trivial, and since $h_n \to
+\infty$ it is a hedgehog. Suppose the fixed point is
accessible from the complement of $K_{\infty}$. Then the lift of
a curve landing at the fixed point gives a curve $\gamma : [0,+\infty) \to W_{\infty}$
such that $\Im \gamma(t) \to +\infty$ as $t \to +\infty$. Let $t_0
> 0$ be such that $\gamma(t) > 10$ for all $t \geq t_0$. Given $T
> t_0$, choose $n \geq 1$ such that $\gamma([t_0, T]) \subset
W_n$, then by the previous Lemma there exists $M$ only depending
on $\gamma(t_0)$ such that $\Im \gamma(T) \leq M$. Hence $\Im
\gamma(t) \leq M$ for all $t \geq t_0$, a contradiction. $\diamond$

\bigskip

\bibliography{hedgehogarea}
\bibliographystyle{alpha}

\medskip

\noindent Ramakrishna Mission Vivekananda University,
Belur Math, WB-711202, India

\noindent email: kingshook@rkmvu.ac.in

\end{document}